\documentclass{amsart}
\usepackage[all]{xy}
\SelectTips{cm}{}

\usepackage{amsmath}
\usepackage{amsthm}
\usepackage{amssymb}

\newtheorem{thm}{Theorem}[section]

\newtheorem{prop}[thm]{Proposition}

\theoremstyle{definition}

\newtheorem{rmk}[thm]{Remark}

\numberwithin{equation}{section}

\def\cal{\mathcal }
\def\cDer{{\cal D}\!er}

\def\cOmega{{\it\Omega}}
\def\cTheta{{\it\Theta}}
\def\cDiff{{\cal D}\!\hbox{\it iff}}

\def\cO{{\cal O}}
\def\cD{{\cal D}}
\def\cE{{\cal E}}

\def\cI{{\cal I}}

\def\cM{{\cal M}}

\def\cP{{\cal P}}
\def\bZ{{\Bbb Z}}
\def\bN{{\Bbb N}}

\def\Gr{{\rm Gr}}
\def\Ch{{\rm Ch}}

\def\Ann{{\rm Ann}}
\def\lMod{{\rm\hbox{-}Mod}}

\def\DR{{\rm DR}}

\def\isom{\cong}
\def\id{{\rm id}}

\title
{Algebraic connections vs. Algebraic {$\cD$}-modules: regularity conditions.}

\author
{Maurizio Cailotto and Luisa Fiorot}
\address{ Dipartimento di Matematica, Universit\`a degli studi di Padova, via Trieste 63, I-35121 Padova, Italy}
\email{maurizio@math.unipd.it, fiorot@math.unipd.it}
\date{}
\begin{document}

\maketitle

\begin{abstract}
This paper is devoted to the comparison of the notions of 
regularity for algebraic connections and (holonomic) regularity 
for algebraic {$\cD$}-modules. 
\end{abstract}

%

\section*{Introduction}
In the dictionary between the language of
(algebraic integrable) connections and that of (algebraic) $\cD$-modules,
the notion of regularity is of great importance, 
and in some sense this justifies different approaches to the definition 
itself. 
In the context of algebraic connections the definition of regularity comes 
from the theory of regular singular points of ordinary differential equations
due to Fuchs: a monic differential operator 
$P=\sum_{i=1}^na_i(x)(x\partial_x)^i$ (with $a_n(x)=1$) is regular at $0$ if 
the coefficients $a_i(x)$ are regular (no poles at $0$), 
or equivalently if in the expression $P=\sum_ib_i(x)\partial_x^i$ (with $b_n(x)=1$)
the coefficients $b_i(x)$ have the property that 
${\rm ord}_0b_i(x)\geq i-n$ ($n$ is the order of $P$). 
In several variables, several notions of regularity (along a polar divisor) 
have been considered. 
The general notion of regularity for an algebraic connection, 
as developed by Manin, Deligne and many other authors, 
is the existence (after suitable localization and completion) 
of a sub-lattice stable under logarithmic derivations. 
In the context of 
$\cD$-modules the notion of regularity, which
generalizes that of regular singular points, is due to
Kashiwara, i.e.: a holonomic $\cD$-module is regular if the 
annihilator of  its graded module w.r.t.~a suitable good filtration is a radical ideal. 
In the ordinary case, that is for analytic functions of one variable, 
these two notions are equivalent by the following elementary argument 
(see \cite{K}). Let $P$ be as before, 
and let us consider the holonomic $\cD$-module $\cM=\cD/\cD P$. 
Then $P$ is regular at $0$ if and only if $\cM$ is regular. 
In fact for the (good) filtration of $\cM$ defined by 
$F_0(\cM)$ being the $\cO$-module generated by 
$u,(x\partial_x)u,\cdots,(x\partial_x)^{n-1}u$, 
and $F_k(\cM)=F_k(\cD)F_0(\cM)$,  we have that 
$x\partial_x$ belongs (and then generates) the annihilator 
of the graded module if and only if 
$x\partial_x F_k(\cM)\subseteq F_k(\cM)$, 
if and only if the operator $P$ has coefficients $a_i(x)$ which are regular. 

In this paper we prove that these two definitions in the general case, 
under suitable conditions, correspond to 
each other in the dictionary. 
Thus we answer a question addressed to us by Andr\'e and Baldassarri 
(as a complement of their book \cite{AB}). 
Even if some authors consider these two notions as equivalent, 
there seems to be no proof of this statement in the literature. 
Hence, this work provides a sequel of our paper \cite{CF}
in the general problem of comparing various notions 
for algebraic connections and for algebraic $\cD$-modules.

\section{Generalities on connections and {$\cD$}-modules}

Let $X$ be a smooth $K$-variety of pure dimension $d_X=\dim X$, 
where $K$ is a field of characteristic $0$.
Following the terminology of \cite[$I\!V$,\S 16]{EGA}, 
we denote by $\cOmega^1_X$ the ${\cO}_X$-module of differentials,
by $\cP_X^1$ the ${\cO}_X$-algebra of principal parts of order one ($1$-jets):
its two structures as ${\cO}_X$-algebra 
(induced by the projections $p_1$, $p_2$ on $X\times X$) 
will be referred to as the 
``left'' and ``right'' structures, and  tensor products will be specified by the
position of the $\cP_X^1$ factor. 
Let us recall that the difference of the inclusions $i_1$, $i_2$ of $\cO_X$ in $\cP_X^1$ 
induced by $p_1$, $p_2$ gives the differential $d=i_2-i_1:\cO_X\to\cOmega^1_X$ 
(i.e. $d(x)=1\otimes x-x\otimes 1$). 
\endgraf 
We also use $\cDer_X$ or $\cTheta_X$ to denote the ${\cO}_X$-module of derivations
(${\cO}_X$-dual of $\cOmega^1_X$, endowed with the usual structure of
Lie-algebra),
and  ${\cD}_X$ to indicate the graded (left) ${\cO}_X$-algebra of differential
operators.
On ${\cD}_X$ we consider the increasing filtration $F$ defined by
the order of differential operators.
Then the associated graded ${\cO}_X$-algebra, denoted by $\Gr\cD_X$,
is commutative and it is generated (as ${\cO}_X$-algebra) by
$\cDer_X\subseteq F^1{\cD}_X$.
\endgraf
For any ${\cO}_X$-module ${\cE}$ we will use the standard notation
$\cP_X^1({\cE})$ for $\cP_X^1\otimes_{{\cO}_X}{\cE}$,
where the tensor product involves the right ${\cO}_X$-module structure
of $\cP_X^1$, while  the ${\cO}_X$-module structure is
given by the left  ${\cO}_X$-module structure  on $\cP_X^1$.

\subsection*{Connections and $\cD$-modules} 
Let $\cE$ be an ${\cO}_X$-module. 
The following supplementary structures on $\cE$ are equivalent: 
\begin{enumerate}
\item
a connection, that is a morphism of abelian sheaves
$\nabla:{\cE}\to\cOmega^1_X\otimes_{{\cO}_X}{\cE}$
which satisfies the Leibniz rule with respect to 
sections of ${\cO}_X$, 
plus the integrability condition, that is $\nabla^2=0$
for the natural extension of $\nabla$ to the De Rham sequence;
\item
an $\cO_X$-linear section $\delta:{\cE}\to\cP^1_X\otimes_{{\cO}_X}{\cE}$
of the canonical morphism $\pi:\cP^1_X\otimes_{{\cO}_X}{\cE}\to\cE$ 
extending to a stratification in the sense of \cite[2.10]{BO};
\item
an  $\cO_X$-linear Lie-algebra homomorphism
$\Delta:\cDer_X\to\cDiff_{X}({\cE})$ (for the usual Lie-algebra
structures), where $\cDiff_{X}({\cE})$ is the sheaf of differential
operators of $\cE$;
\item
a structure of left  
$\cD_X$-module on $\cE$.
\end{enumerate}
\endgraf
The dictionary between these equivalent structures is
well explained in \cite[2.9, 2.11, 2.15]{BO}; let us give a sketch.
If $c=c_X({\cE}):{\cE}\to\cP^1_X\otimes_{{\cO}_X}{\cE}$ denotes the
inclusion induced by $i_2$ ($1$-jets), 
then $\delta=c-\nabla$ and $\nabla=c-\delta$.
For  any $\partial$ section of $\cDer_X$
the morphism $\Delta$ is defined by
$\Delta_{\partial}=(\partial\otimes\id)\circ\nabla$,
i.e. $\Delta_{\partial}(e)=\langle\partial,\nabla(e)\rangle$.
On the other hand,  the reconstruction of $\nabla$ from $\Delta$ involves a
description
using local coordinates $x_i$ on $X$ 
($dx_i$ and $\partial_i$ are the dual bases of differentials and derivations):
if $e$ is a section of $\cE$, then
$\nabla(e)=\sum_idx_i\otimes\Delta_{\partial_i}(e)$.
\endgraf
The morphism $\Delta$ is equivalent to the data of a
left $\cD_X$-module structure on  $\cE$ since it extends
to a left action of $\cD_X$ on $\cE$ (see \cite[VI,1.6]{Bo}).
In fact the datum of a connection (without the integrability condition) 
is equivalent to the datum of a section of $\pi$ 
(without further conditions), as explained in \cite[I,2.3]{D}, 
and in that correspondence, since $K$ is of characteristic $0$, 
integrable connections correspond  to sections extending to stratifications 
(see \cite[2.15]{BO}). 
From now on, the word connection means integrable connection, 
that is connection satisfying the integrability condition. 

\subsection*{Morphisms} 
A morphism of connections on $X$ is an ${\cO}_X$-linear morphism
$h:{\cE}\to{\cE}'$ compatible with the data, that is, such that
$\nabla'\circ h=({\rm id}\otimes h)\circ\nabla$,
or $\delta'\circ h=({\rm id}\otimes h)\circ\delta$,
or equivalently
$\Delta'_\partial\circ h=h\circ\Delta_\partial$
for any section $\partial$ of $\cDer_X$,
or finally which is $\cD_X$-linear.

\subsection*{Coherence and quasi-coherence conditions} 

The connection $\cE$ is said to be quasi-coherent (resp. coherent)
if $\cE$ enjoys the corresponding property as  ${\cO}_X$-module.
Recall that coherence implies locally freeness for integrable connections
(see \cite[2.17]{BO}).

Let denote by ${\rm MIC}(X)$ (resp.${\rm MIC}_{qc}(X)$, resp. ${\rm MIC}_{c}(X)$)
the category of integrable
(resp. quasi-coherent, resp. coherent so locally free of finite type) 
connections.

For us a ${\cD}_X$-module is a left algebraic ${\cD}_X$-module
and we denote this category by $\cD_X\lMod$.
Let $i:\cO_X\hookrightarrow \cD_X$ be the usual inclusion
of $\cO_X$ into $\cD_X$.

A $\cD_X$-module $\cM$ is coherent if 
for any $x\in X$ there exists an affine neighborhood $U$
and an exact sequence 
$$\xymatrix{
\cD_X(U)^q \ar[r] & \cD_X(U)^p  \ar[r] &  \cM(U)   \ar[r] &  0}
$$
(see \cite[VI,1.4]{Bo}).
We denote by $\cD_X\lMod_{c}$ the category of coherent
$\cD_X$-modules (warning: they may not be coherent as $\cO_X$-modules,
for example $\cD_X$ is coherent as $\cD_X$-module but it is only
quasi-coherent as $\cO_X$-module).
Any coherent $\cD_X$-module is
quasi-coherent as ${\cO}_X$-module (see \cite[VI.2.11]{Bo}).
Moreover a $\cD_X$-module which is coherent as $\cO_X$-module is locally
$\cO_X$-free of finite type (\cite[VI.1.7]{Bo}) and we denote by
$\cD_X\lMod_{\cO_X\hbox{-}c}$
the full subcategory of $\cD_X\lMod$ whose objects are 
$\cO_X$-coherent $\cD_X$-modules.
\endgraf

A $\cD_X$-module $\cM$ is quasi-coherent
if for any $x\in X$ there exists an affine neighborhood $U$
and an exact sequence 
$$\xymatrix{
\cD_X(U)^{(I)}  \ar[r] & \cD_X(U)^{(J)} \ar[r] & \cM(U)  \ar[r] &  0}
$$
where $I, J$ are arbitrary set of indexes and $\cD(U)^{(I)}$ represents
the direct sum of $\cD(U)$ indexed by $I$.
Let $\cD_X\lMod_{qc}$ be the category of quasi-coherent
${\cD}_X$-modules.
Any quasi-coherent
${\cD}_X$-module is quasi-coherent as $\cO_X$-module
(because $\cD_X$ is a quasi-coherent $\cO_X$-module
and direct sums of quasi-coherent $\cO_X$-modules
are quasi-coherent $\cO_X$-modules).
Moreover any $\cD_X$-module which is quasi-coherent as $\cO_X$-module
is also quasi-coherent as $\cD_X$-module.
In fact for any $x\in X$ there exists an affine  neighborhood $U$
and an epimorphism  
$$\xymatrix{
\cO_X(U)^{(J)}\ar[r]^{g} &\cM(U)  \ar[r] &0.
}
$$
Let $\overline{g}$ be the morphism obtained by extension of scalars
from $\cO_X(U)$ to $\cD_X(U)$.
Then   $\overline{g}$ too is an epimorphism whose kernel will be denoted by $K$.
This is a $\cD_X(U)$-module and there exist $I$ and $f$ such that the morphism
$f:\cO_X(U)^{(I)}\longrightarrow K$
is surjective.
As before let $\overline{f}$ be the morphism obtained from $f$ 
extending the scalars to $\cD_X(U)$.
We obtain an exact sequence
$$
\xymatrix{
\cD_X(U)^{(I)}  \ar[r]& \cD_X(U)^{(J)} \ar[r] & \cM(U)  \ar[r] &  0}
$$
which proves that $\cM$ is a quasi-coherent  $\cD_X$-module.
 
The full subcategory of ${\cD}_X$-modules which are
quasi-coherent (resp. coherent) as $\cO_X$-modules
is isomorphic to the category of quasi-coherent (resp. coherent)
connections.
So we have the following commutative diagram whose
horizontal arrows are isomorphisms of categories:
$$
\xymatrix@R=15pt{
{\rm MIC}_{c}(X) \ar[r]\ar[d] & \cD_X\lMod_{\cO_X\hbox{-}c} \ar[d] \\
{\rm MIC}_{qc}(X) \ar[r]\ar[d] & \cD_X\lMod_{qc}\ar[d] \\
{\rm MIC}(X) \ar[r] &  \cD_X\lMod. \\
}
$$
In the following we consider only
quasi-coherent $\cD_X$-modules.

\section{Definitions of regularity}

\subsection*{Good filtrations of ${\cD}_X$-modules}\label{ssGF} 
A filtration $F^i({\cM})$ of a ${\cD}_X$-module
$\cM$ is an increasing $\bZ$-indexed
family of coherent sub-${\cO}_X$-modules of $\cM$ such that
$F^i{\cM}=0$ for $i\ll 0$, $\cM$ is the union of all the $F^i{\cM}$ and
$F^i{\cD}_X\, F^j{\cM}\subseteq F^{i{+}j}{\cM}$.
The filtration is said to be good (or coherent) if one of the following
equivalent conditions holds: 
\begin{description}
\item[{\bf 
(i)}]
for $j\gg 0$ and all $i\in\bN$
we have $F^i{\cD}_X\, F^j{\cM}=F^{i{+}j}{\cM}$;
\item[{\bf 
(ii)}]
the associated graded module $\Gr_F{\cM}=\bigoplus_{i\in\bZ}\Gr^i_F{\cM}$
(where $\Gr^i_F{\cM}=F^i{\cM}/F^{i{-}1}{\cM}$) is a
coherent $\Gr{\cD}_X$-module.
\end{description}
We recall that in the algebraic setting (and unlike the analytic case)
any coherent  ${\cD}_X$-module admits a {\it global} good filtration
(\cite[I.2.5.4]{Me}).

\subsection*{Characteristic variety of ${\cD}_X$-modules}
Let $T^\ast X={\bf V}((\cOmega^1_X)^\vee)$ be the cotangent bundle of $X$
(we use in general the terminology of \cite[II]{EGA}).
We denote by $\pi=\pi_X$ the canonical morphism of $K$-varieties
$T^\ast X\to X$ and by $\iota=\iota_X:X\to T^\ast X$
the zero section of $\pi$, whose image is $T^\ast_XX$.
\endgraf
For any ${\cD}_X$-module $\cM$ and any good filtration $F$ on it,
the graded module $\Gr_F{\cM}$ is a
${\cO}_{T^\ast X}=\Gr{\cD}_X$-module.
The characteristic variety 
$\Ch\cM$ of $\cM$ is defined as the support in $T^\ast X$ of $\Gr_F{\cM}$,
that is the closed subset of $T^\ast X$ corresponding to the
annihilator $\cI_F(\cM)={\rm Ann}_{\Gr{\cD}_X}(\Gr_F{\cM})$
of  $\Gr_F{\cM}$ in ${\cO}_{T^\ast X}$.
We recall that the ideal $\cI_F(\cM)$ depends on the filtration $F$,
but the characteristic variety $\Ch\cM$ does not,
that is, the radical of $\cI_F(\cM)$ is independent of $F$
(see for example \cite{Ga} and \cite[2.6]{K}).
Moreover the characteristic variety of a $\cD_X$-module
is always a conical involutive closed subset in $T^\ast X$ (see \cite[VI.1.9]{Bo},
\cite[I.2.3;2.5]{Me}, \cite{Ga}), 
and in particular the Bernstein inequality holds: $\dim\Ch\cM\geq\dim X$ 
(see \cite[VI.1.10]{Bo}, \cite[I.2.3.4;2.5]{Me}).
\endgraf
A ${\cD}_X$-module $\cM$ is $\cO_X$-coherent if and only if
$\Ch\cM=T^\ast_XX$.
\endgraf

\subsection*{Holonomic ${\cD}_X$-modules}
A ${\cD}_X$-module $\cM$ is said to be holonomic if
$\dim\Ch\cM\leq\dim X$ (so that the equality holds, and the
characteristic variety has the minimal possible dimension).
We denote by $\cD_X\lMod_h$ the category of holonomic
$\cD_X$-modules (as a full subcategory of $\cD_X\lMod$).
Since any $\cO_X$-coherent $\cD_X$-module is holonomic,
$\cD_X\lMod_{\cO_X\hbox{-}c}$ is a full subcategory of $\cD_X\lMod_h$.
Notice that a ${\cD}_X$-module is holonomic if and only if 
its characteristic variety is lagrangian 
(and so a union of conormal varieties).

\subsection*{Regularity for holonomic {${\cD}_X$}-modules}

Following Kashiwara (see  \cite[5.2]{K}),
a holonomic $\cD_X$-module $\cM$ is said to be regular, or to
have regular singularities (or to be RS)
if it admits a good filtration  $F$ such that $\cI_F(\cM)$ is a radical
ideal,
or equivalently the (reduced) ideal $\cI(\Ch\cM)$ of $\Ch\cM$ annihilates
$\Gr_F\cM$.

Let $\cM$ be an $\cO_X$-coherent $\cD_X$-module.
 Then $\cM$  belongs to
$\cD_X\lMod_h$ and it always has  regular singularities; in fact
we can take $F^i(\cM)=\cM$ for any $i\geq 0$ and $F^i(\cM)=0$ if $i<0$.
Then $\cI_F(\cM)=\bigoplus_{k\geq 1}\Gr^k\cD_X$
which is a radical ideal.

Let $\cM$ be a $\cD_X$-module.
A point $x\in X$ is called a singularity for $\cM$ if
$\big(\pi^{-1}(x)\smallsetminus 
T^\ast_XX \big)\cap\Ch(\cM)\neq
\emptyset$.
In particular, 
a $\cD_X$-module which is $\cO_X$-coherent  has no singular points.

\subsection*{Regularity for connections} 

Let $X$ be a smooth $K$-variety and let $Z$ be a smooth irreducible
hypersurface of $X$. A connection $(\cE ,\nabla)$ on $U=X\smallsetminus Z$
is said to be regular along $Z$ if (and only if) $E=\cE_{\eta_X}$ 
($\eta_X$ is the generic point of $X$, and $U$) is
a $\kappa(X)/K$-differential module regular
at the divisorial valuation $v$ corresponding to $Z$, 
that is, the completion of $E$ w.r.t. $v$ admits 
a sub-$\widehat{\cO}_{X,\eta_Z}$-lattice 
stable under $x\partial_x$ where $x$ is a local equation for $Z$ 
(a generator for the ideal $\cI_Z$ of $Z$ in $\cO_X$), 
and $\partial_x$ is a derivation transversal to $Z$ 
(i.e. such that $\partial_x(m_{X,Z})\not\subset m_{X,Z}$ 
where $m_{X,Z}=\widehat{\cI}_{Z,\eta_Z}$) 
satisfying $\partial_x(x)=1$. 
\endgraf
Let $X$ be a smooth $K$-variety.
A connection $(\cE,\nabla)$ on $X$ is said to be regular if
$E=\cE_{\eta_X}$ is a $\kappa(X)/K$-differential module regular
at any divisorial valuation of $\kappa(X)/K$. 
\endgraf 
This definition shows immediately that the notion of 
regularity is a birational invariant. 
It is useful to have a more concrete characterization: 
$(\cE,\nabla)$ on $X$ is regular if there exists a normal compactification 
$\overline{X}$ of $X$ such that the connection is regular along any 
component of the boundary $Z=\overline{X}\setminus X$ which is 
of codimension one in $\overline{X}$. 
This characterization is easier to prove if we suppose $Z$ to be a normal crossing divisor 
(in that case a proof can be done with a logarithmic differential argument); 
while the general case has been proved with analytic methods by Deligne 
(in \cite{D} the proof of this criterion contains a mistake, 
and a correct proof is given in the ``erratum'' of 1971), 
then with algebraic methods by Andr\'e (\cite{A}). 
\endgraf
Whenever $Z=\overline{X}\setminus X$ is 
a normal crossing divisor, a connection $(\cE ,\nabla)$ is regular 
if and only if there exists an extension $\widetilde\cE$ of $\cE$ to $X$, 
and $\widetilde\nabla$ of $\nabla$ with logarithmic poles along $Z$.  
Such an extension is unique if the eigenvalues of the residues of 
the connection are forced to belong to the image of a section 
$\tau$ of the canonical projection $K\to K/\bZ$ 
($\tau$-extension of Deligne, 
constructed in \cite{D} with analytic methods, 
then in \cite{AB} with algebraic methods).

\subsection*{Connections with poles}
Let  $X$ be a smooth $K$-variety, $Z$ a
divisor with 
normal crossings in $X$
(we denote by $j$ the inclusion of $U$ in $X$)
and $\cE$ an $\cO_U$-coherent $\cD_U$-module
(so that it is locally free of finite rank as $\cO_U$-module).
Let $\overline{\cE}$ be a coherent $\cO_X$-module contained in $j_\ast(\cE)$
such that
$j^{-1}(\overline{\cE})=\cE$.
We call such an $\overline{\cE}$ a coherent extension of $\cE$ to $X$.
We have 
$j_\ast\cE\isom\overline{\cE}(\ast Z):=
\lim\limits_{\longrightarrow}{\!}_i\,\cI_Z^{-i}\overline{\cE}$
for any coherent extension $\overline{\cE}$ of $\cE$ as before and
in particular
$j_\ast\cO_U\isom\cO_X(\ast Z):=
\lim\limits_{\longrightarrow}{\!}_i\,\cI_Z^{-i}\cO_X$
and
$j_\ast\cOmega^1_U\isom\cOmega^1_X(\ast Z):=
\lim\limits_{\longrightarrow}{\!}_i\,\cI_Z^{-i}\cOmega^1_X$.
\endgraf
 
Let us denote by
$\cTheta_{X,Z}\subset \cTheta_X$
the sheaf of derivations
which respect the ideal $\cI_Z$ 
and  by $\cD_{X,Z}$ the
 sub-$\cO_X$-algebra of $\cD_{X}$
generated by the derivations $\cTheta_{X,Z}$.

Let  $j:U\hookrightarrow X$ be an open immersion, and suppose 
 that $Z=X\smallsetminus U$ is a divisor with normal crossings. 
 Let consider the sheaf $\cO_U$ endowed with the trivial connection. 
 Then $j_\ast(\cO_U)=\cO_X(\ast Z)$ is a regular holonomic $\cD_X$-module 
 and $\Ch(j_\ast(\cO_U))=V(\cTheta_{X,Z}\Gr(\cD_X))$, 
 where $\cTheta_{X,Z}\Gr(\cD_X)$ is the ideal generated 
 by $\cTheta_{X,Z}$ in $\Gr(\cD_X)$. 
 More precisely, let $F$ be the good filtration on $j_\ast(\cO_U)$ 
 which is zero for negative degrees and is generated 
(as $\cO_X$-module) in degree $i$ by 
 sections of $j_\ast(\cO_U)$ with poles of order $i$ 
 on $Z$. 
 Then $\Ann\Gr_F(j_\ast(\cO_U))=\cTheta_{X,Z}\Gr(\cD_X)$. 
 \endgraf 
 In fact, let us suppose $Z$ has locally equation $x_1\cdots x_r$ 
 using local coordinates $x_1,\dots,x_n$ in $X$.  
 Then $\cD_{X,i}F^j(j_\ast(\cO_U))=F^{j+i}(j_\ast(\cO_U))$ 
 for any $j\geq r$, so that $F$ is a good filtration. 
 Clearly $\cTheta_{X,Z}\Gr(\cD_X)$ 
 is contained in $\Ann\Gr_F(j_\ast(\cO_U))$. 
On the other side, a local computation shows immediately that 
any section $s$ of $\Ann\Gr_F(j_\ast(\cO_U))$ 
belongs to $\cTheta_{X,Z}\Gr(\cD_X)$ 
(for example, applying $s$ to 
$1\over x_i$ for all $i=1,\cdots,r$). 

\begin{prop}
Let $(\cE,\nabla)$ be a coherent connection on $U$.
For any good filtration $F$ on $j_{\ast}(\cE)$,  
we have 
$\Ann\Gr_F(j_\ast(\cE))\subseteq \cTheta_{X,Z}\Gr(\cD_X)$.
\end{prop}
%
%
{\bf Proof.}
Since $j_\ast(\cE)$ is a holonomic $\cD_X$-module, 
its characteristic variety is a conical lagrangian subvariety of $T^\ast X$. 
Since it has poles along any component of $Z$, 
the characteristic variety contains $T^\ast_ZX$. 
Therefore, for any good filtration $F$, we have 
$\Ann\Gr_F(j_\ast(\cE))
\subseteq \sqrt{\Ann\Gr_F(j_\ast(\cE))}
\subseteq \cTheta_{X,Z}\Gr(\cD_X)$. 
\hfill$\square$

\begin{rmk}
In the previous proof we have used  that $\Ch(j_\ast\cE)\supseteq T^\ast_ZX$.
This  fact is hard to prove directly via
local computations. 
Following a remark of C.~Sabbah we may obtain a better understanding 
of the situation: using the exact sequence of perverse sheaves 
$$ 
0\longrightarrow 
{\rm Irr}_Z(j_\ast\cE)\longrightarrow 
\DR(j_\ast\cE)\longrightarrow 
Rj_\ast\DR(\cE)\longrightarrow 0
$$ 
due to Mebkhout (see \cite{Me2}) we see that 
the characteristic variety of the middle term $\Ch(j_\ast\cE)=\Ch(\DR(j_\ast\cE))$ 
is the sum of $\Ch(Rj_\ast\DR(\cE))$ 
(which is by Riemann-Hilbert the characteristic variety of a regular 
holonomic $\cD$-module, so that it is exactly $T^\ast_ZX$), 
and $\Ch({\rm Irr}_Z(j_\ast\cE))$ 
(and the irregularity sheaf need not to be adapted to the natural 
stratification of $Z$). 
\endgraf 
For example if we consider the irregular modules with solution $e^{x/y}$ 
(using $x,y$ local coordinates of the plane, it has poles only along $Z$ defined by $y=0$), 
then the characteristic variety has three 
components: $T^\ast_XX$, $T^\ast_ZX$ and $T^\ast_0X$ 
(due to the lack of a good formal structure at $0$). 
\end{rmk}

\section{Comparison}

We now compare the notion of regularity for connections and $\cD$-modules.
Let us remark that the notion of regular connection  $(\cE,\nabla)$ (with
$\cE$ a coherent $\cO_U$-module) takes in account 
the so called regularity at infinity  where the connection has poles.
If we consider $(\cE,\nabla)$ as a $\cD_U$-module it is always regular in the
sense of Kashiwara (as previously noticed). Hence we need to pass to a compactification in order to compare correctly these notions.

\begin{thm}\label{compT}
Let $U$ be a smooth $K$-variety and let $j:U\hookrightarrow X$ be an open
dense immersion
where $X$ is a smooth proper $K$-variety and $Z:=X\smallsetminus U$ 
is a divisor with strict normal crossings.
Let $(\cE,\nabla)$ be a coherent connection on $U$.
Then the following are equivalent:
\begin{enumerate}
\item $(\cE,\nabla)$ is regular;
\item $j_\ast\cE$ is a regular holonomic $\cD_X$-module.
\end{enumerate}
\end{thm}
{\bf Proof.}

{\bf(1)$\Rightarrow$\bf(2)}.
Let $(\cE,\nabla)$ be a regular connection on $U$ (along $Z$), 
and consider a $\tau$-extension $(\widetilde\cE,\widetilde\nabla)$ to $X$ 
with logarithmic poles along $Z$. 
We have $j_\ast\cE=\widetilde\cE\otimes_{\cO_X}\cO_X(\ast Z)$, 
and we define on it the filtration given by 
$F^0(j_\ast\cE)=\widetilde\cE\otimes_{\cO_X}\cO_X(dZ)$, 
where $d$ is the number of components of $Z$, 
and $F^i(j_\ast\cE)=\cD_{X,i}F^0(j_\ast\cE)$. 
It gives 
a good filtration, 
and the annihilator of the associated graded module 
contains the ideal which preserves $\widetilde\cE$, 
so that it contains the whole ideal $\cTheta_{X,Z}\Gr(\cD_X)$. 
Since it cannot be bigger by the previous discussion, the annihilator is just 
$\cTheta_{X,Z}\Gr(\cD_X)$, which is a radical ideal. 

{\bf(2)$\Rightarrow$\bf(1)}. 
Let $(\cE,\nabla)$ be a connection on $U$, 
and suppose that $j_\ast\cE$ is a regular holonomic $\cD_X$-module. 
We have to prove that the connection is regular along any 
one codimensional component $Z_i$ of $Z$. 
By hypothesis there exists a good filtration $F^i$ on $j_\ast\cE$ 
with the property that the annihilator of $\Gr_F(j_\ast\cE)$ is a radical ideal 
of $\Gr(\cD_X)$.
Up to a shift on the filtration we may suppose $\widetilde\cE=F^0(j_\ast\cE)$, $F^i(j_\ast\cE)=\cD_{X,i}\widetilde\cE$ for $i\geq 0$ and $F^i(j_\ast\cE)=0$ for $i<0$.
Now, $j_\ast\cE=\widetilde\cE\otimes_{\cO_X}\cO_X(\ast Z)$, 
and let $x_1,\dots,x_n$ be local coordinates 
such that 
$x_1\cdots x_d$ is a local equation for $Z$.
Hence we know that $\partial_{x_i}$ acts on a trivialization of $\widetilde\cE_{\eta_{Z_i}}$ 
via a matrix with poles in $x_1\cdots x_d$. 
Let us denote by
$s$ the maximal order of these poles. 
Therefore 
$(x_1\cdots x_d)^s\partial_{x_i}$ belongs to the annihilator 
of $\Gr_F(j_\ast\cE)$ since for any $i=1,\dots ,d$ we have
$(x_1\cdots x_d)^s\partial_{x_i}\widetilde\cE_{\eta_{Z_i}}\subseteq
\widetilde\cE_{\eta_{Z_i}}$ and so 
$(x_1\cdots x_d)^s\partial_{x_i}\cD_{X,k}\widetilde\cE_{\eta_{Z_i}}\subseteq \cD_{X,k}\widetilde\cE_{\eta_{Z_i}}$ which proves that 
$(x_1\cdots x_d)^s\partial_{x_i}\Gr_F(j_\ast\cE)=0$.
Now the radicality of the  annihilator implies that also 
$x_1\cdots x_d\partial_{x_i}$ belongs to the annihilator. 
In particular $\widetilde\cE_{\eta_{Z_i}}$ is stable under $x_i\partial_{x_i}$. 
Taking the completion w.r.t. the valuation induced by $Z_i$ 
we have an $\widehat{\cO}_{{X},\eta_{Z_i}}$-lattice 
stable under $x_i\partial_{x_i}$ as required. 
\hfill$\square$

\begin{rmk}
We note that the proof of $(2)\Rightarrow (1)$ in \ref{compT} generalizes to
the case of
$Z$ a general hypersurface
simply restricting to $X\setminus{\rm Sing}(Z)$. 
By contrast, the implication $(1)\Rightarrow (2)$ uses 
essentially the existence of $\tau$-extensions, 
which requires to have a normal crossing divisor. 
\end{rmk}

\begin{rmk}
We may try to prove the theorem by reduction to the case 
of dimension one 
(i.e. for curves: in that case the equivalence of the definitions is sketched in \cite{K}
and in the introduction of this paper), 
but it seem to be difficult to prove that 
the Kashiwara definition of regular holonomic $\cD_X$-modules 
can be recovered in terms of curves. 
\end{rmk}

\begin{rmk}
The definition of regular singularity used 
in \cite{KT} (in the general microlocal context) 
or \cite{Be} (in the algebraic $\cD$-modules context) 
is clearly equivalent to the notion of regularity of the 
correspondent object of ${\rm MIC}$.  
Hence, in its $\cD$-module counterparts,  
it corresponds to the regularity of its direct image 
by an open immersion to a proper variety with 
complement being a normal crossing divisor. 
\end{rmk}

\begin{prop}
In the above situation
the following conditions are equivalent:
\item{\bf(a)}
the $\cD_X$-module $j_\ast\cE$ has regular singularities;
\item{\bf(b)}
there exists an $\cO_X$-coherent extension
$\overline\cE$ of $\cE$ to $X$ which is a $\cD_{X,Z}$-module;
\item{\bf(c)}
there exists an $\cO_X$-coherent  extension
$\overline\cE$ of $\cE$ to $X$ such that
${\rm Im}(\cD_{X,Z}\times \overline\cE\rightarrow j_{\ast}(\cE))$
is $\cO_X$-coherent;
\item{\bf(d)}
for any $\cO_X$-coherent extension
$\overline\cE$ of $\cE$ to $X$ the $\cO_X$-module
${\rm Im}(\cD_{X,Z}\times \overline\cE\rightarrow j_{\ast}(\cE))$
is $\cO_X$-coherent.
\endgraf
\end{prop}
 
It is the analog of \cite[I,3.3.4]{AB}.


\end{document}